\documentclass[a4paper]{amsart}

\RequirePackage{amsmath}
\RequirePackage{bm}
\RequirePackage{amssymb}
\RequirePackage{upref}
\RequirePackage{amsthm}
\RequirePackage{enumerate}
\RequirePackage{pb-diagram}
\RequirePackage{amsfonts}
\RequirePackage[mathscr]{eucal}
\RequirePackage{verbatim}
\RequirePackage{xr}
\RequirePackage{graphicx}
\RequirePackage{floatflt}
\usepackage{calc}
\usepackage{xspace}

\newcommand{\cf}{cf.\@\xspace}

\newcommand{\al}{\alpha}
\newcommand{\bet}{\beta}

\newcommand{\de}{\delta }
\newcommand{\e}{\epsilon}

\newcommand{\f}{\varphi}
\newcommand{\h}{\eta}

\newcommand{\ka}{\kappa}
\newcommand{\lam}{\lambda}

\newcommand{\n}{\nu}

\newcommand{\s}{\sigma}
\newcommand{\x}{\xi}

\newcommand{\C}{\varGamma}

\newcommand{\F}{\varPhi}

\newcommand{\Om}{\varOmega}

\newcommand{\nbdd}{\nobreakdash--}

\newcommand{\fv}[2]{#1\hspace{0pt}_{|_{#2}}}

\newcommand{\so}{{\mc S_0}}

\newcommand{\const}{\tup{const}}
\newcommand{\ndash}{\nobreakdash--}

\newcommand{\msp[1]}[1]{\mspace{#1mu}}

\newcommand{\R}[1][n+1]{{\protect\mathbb R}^{#1}}

\newcommand{\N}{{\protect\mathbb N}}

\newcommand{\eR}{\stackrel{\lower1ex \hbox{\rule{6.5pt}{0.5pt}}}{\msp[3]\R[]}}
\newcommand{\eN}{\stackrel{\lower1ex \hbox{\rule{6.5pt}{0.5pt}}}{\msp[1]\N}}
\newcommand{\eO}{\stackrel{\lower1ex
\hbox{\rule{6pt}{0.5pt}}}{\msc O}}

\DeclareMathOperator{\graph}{graph}

\newcommand\pa{\partial}
\newcommand\pde[2]{\frac {\partial#1}{\partial#2}}
 
\newcommand{\A}{\forall}

\newcommand{\set}[2]{\{\,#1\colon #2\,\}}
\newcommand{\uu}{\cup}
\newcommand{\ii}{\cap}
\newcommand{\uuu}{\bigcup}

\newcommand{\uud}{ \stackrel{\lower 1ex \hbox {.}}{\uu}}
\newcommand{\uuud}[1]{ \stackrel{\lower 1ex \hbox {.}}{\uuu_{#1}}}
\newcommand\su{\subset}

\newcommand{\sminus}[1][28]{\raise 0.#1ex\hbox{$\scriptstyle\setminus$}}

\newcommand{\abs}[1]{\lvert#1\rvert}

\newcommand{\norm}[1]{\lVert#1\rVert}

\newcommand{\spd}[2]{\protect\langle #1,#2\protect\rangle}

\newcommand{\riema}[4]{{\bar R}_{#1#2#3#4}}

\newcommand{\tit}{\textit}

\newcommand{\tup}{\textup}

\newcommand{\mc}{\protect\mathcal}
\newcommand{\msc}{\protect\mathscr}

\providecommand{\bysame}{\makebox[3em]{\hrulefill}\thinspace}

\newcommand{\bt}{\begin{thm}}
\newcommand{\bl}{\begin{lem}}
\newcommand{\bc}{\begin{cor}}
\newcommand{\bd}{\begin{definition}}
\newcommand{\bpp}{\begin{prop}}
\newcommand{\br}{\begin{rem}}
\newcommand{\bn}{\begin{note}}
\newcommand{\be}{\begin{ex}}
\newcommand{\bes}{\begin{exs}}
\newcommand{\bb}{\begin{example}}
\newcommand{\bbs}{\begin{examples}}
\newcommand{\ba}{\begin{axiom}}
\newcommand{\bas}{\begin{assumption}}

\newcommand{\et}{\end{thm}}
\newcommand{\el}{\end{lem}}
\newcommand{\ec}{\end{cor}}
\newcommand{\ed}{\end{definition}}
\newcommand{\epp}{\end{prop}}
\newcommand{\er}{\end{rem}}
\newcommand{\en}{\end{note}}
\newcommand{\ee}{\end{ex}}
\newcommand{\ees}{\end{exs}}
\newcommand{\eb}{\end{example}}
\newcommand{\ebs}{\end{examples}}
\newcommand{\ea}{\end{axiom}}
\newcommand{\eas}{\end{assumption}}

\newcommand{\bp}{\begin{proof}}
\newcommand{\ep}{\end{proof}}
\newcommand{\eps}{\renewcommand{\qed}{}\end{proof}}

\newcommand{\bal}{\begin{align}}

\newcommand{\bi}[1][1.]{\begin{enumerate}[\upshape #1]}
\newcommand{\bia}[1][(1)]{\begin{enumerate}[\upshape #1]}
\newcommand{\bin}[1][1]{\begin{enumerate}[\upshape\bfseries #1]}
\newcommand{\bir}[1][(i)]{\begin{enumerate}[\upshape #1]}
\newcommand{\bic}[1][(i)]{\begin{enumerate}[\upshape\hspace{2\cma}#1]}
\newcommand{\bis}[2][1.]{\begin{enumerate}[\upshape\hspace{#2\parindent}#1]}
\newcommand{\ei}{\end{enumerate}}

\newcommand\ndots{\raise 0.47ex \hbox {,}\hskip0.06em\cdots %
     \raise 0.47ex \hbox {,}\hskip0.06em} 

\newcommand{\q}{\quad}
\newcommand{\qq}{\qquad}

\newcommand{\hp}{\hphantom}

\newcommand\nd{\noindent}

\newskip\Csmallskipamount                                                
\Csmallskipamount=\smallskipamount
\newskip\Cmedskipamount
\Cmedskipamount=\medskipamount
\newskip\Cbigskipamount
\Cbigskipamount=\bigskipamount

\newcommand\cvs{\vspace\Csmallskipamount}   
\newcommand\cvm{\vspace\Cmedskipamount}

\newskip\csa
\csa=\smallskipamount

\newskip\cma
\cma=\medskipamount

\newskip\cba
\cba=\bigskipamount

\newdimen\spt
\newcommand\citem{\cvs\advance\itemno by
1{(\romannumeral\the\itemno})\hskip3pt}
\newcommand{\bitem}{\cvm\nd\advance\itemno by
1{\bf\the\itemno}\hspace{\cma}}

\newcount\itemno
\newcommand{\lae}[1]{\label{E:#1}}
\newcommand{\lat}[1]{\label{T:#1}}
\newcommand{\lal}[1]{\label{L:#1}}

\newcommand{\lar}[1]{\label{R:#1}}

\newcommand{\laas}[1]{\label{Ass:#1}}

\newcommand{\rt}[1]{Theorem~\ref{T:#1}}

\newcommand{\ras}[1]{Assumption~\ref{Ass:#1}}

\newcommand{\re}[1]{\eqref{E:#1}}

\newskip\thmskip
\thmskip=\parindent

\newskip\hsk
\newenvironment{hinw}{\labelsep=0pt\begin{list}{}{\labelsep=0pt\itemindent=0pt\labelwidth=0pt\leftmargin=\parindent\rightmargin=0pt\partopsep=\cba}%
\item\it\nopagebreak\nopagebreak}%
{\end{list}}

\newcommand\bh{\begin{hinw}}
\newcommand{\eh}{\end{hinw}}

\newtheoremstyle{normal}
  {\cba}
  {\cba}
  {}
  {\thmskip}
  {\bfseries}
  {.}
  {\hsk}
  {}

\newtheoremstyle{abschnitt}
  {\cba}
  {\cba}
  {}
  {\thmskip}
  {\bfseries}
  {.}
  {\hsk}
  {}

\newtheoremstyle{italic}
  {\cba}
  {\cba}
  {\itshape}
  {\thmskip}
  {\bfseries}
  {.}
  {\hsk}
  {}

\newtheoremstyle{aufgaben}
  {\cba}
  {\cba}
  {}
  {}
  {\normalsize\bfseries}
  {.}
  {\hsk}
  {}

\newtheoremstyle{break}
  {\cba}
  {\cba}
  {\itshape}
  {}
  {\bfseries}
  {.}
  {\newline}
  {}

\swapnumbers
\theoremstyle{italic}
\newtheorem{thm}[subsection]{Theorem}
\newtheorem{lem}[subsection]{Lemma}
\newtheorem{prop}[subsection]{Proposition}
\newtheorem{cor}[subsection]{Corollary}

\theoremstyle{normal}
\newtheorem{rem}[subsection]{Remark}
\newtheorem{definition}[subsection]{Definition}
\newtheorem{example}[subsection]{Example}
\newtheorem{examples}[subsection]{Examples}
\newtheorem{ex}[subsection]{Exercise}
\newtheorem{note}[subsection]{}
\newtheorem{axiom}[subsection]{Axiom}
\newtheorem{assumption}[subsection]{Assumption}

\theoremstyle{aufgaben}
\newtheorem{exs}[subsection]{Exercises}

\swapnumbers

\numberwithin{equation}{section}
\numberwithin{figure}{section}

\newenvironment{textequation}[1][0.8]
{\begin{equation}
\begin{aligned}
\begin{minipage}{#1\linewidth}}
{\end{minipage}
\end{aligned}
\end{equation}
\ignorespacesafterend}

\newcommand{\btext}{\begin{textequation}}
\newcommand{\etext}{\end{textequation}}

\usepackage[german,english]{babel}
\usepackage{graphicx}
\RequirePackage{amsmath}
\RequirePackage{bm}
\RequirePackage{amssymb}
\RequirePackage{upref}
\RequirePackage{amsthm}
\RequirePackage{enumerate}
\RequirePackage{pb-diagram}
\RequirePackage{amsfonts}
\RequirePackage[mathscr]{eucal}
\RequirePackage{verbatim}
\RequirePackage{xr}
\RequirePackage{graphicx}
\usepackage{calc}
\usepackage{xspace}

\makeatletter
\RequirePackage{color}
\newcommand{\ann}[1]{\renewcommand{\@makefnmark}{\mbox{$^{\color{red}{\@thefnmark}}$}}%
\footnote {#1}}
\makeatother

\RequirePackage{upref}
\RequirePackage{amsthm}
\RequirePackage{enumerate}
\usepackage[mathscr]{eucal}

\usepackage{xr-hyper}

\usepackage{calc}

\newlength{\oddsidemarginlength}
\newlength{\topmarginlength}

\newcounter{numberoflines}
\newcounter{tempcc}
\oddsidemargin=\oddsidemarginlength
\evensidemargin=\oddsidemargin
\topmargin=\topmarginlength
\begin{document}

\flushbottom


\title[Curvature estimates]{Curvature estimates for Weingarten hypersurfaces in Riemannian manifolds}

\author{Claus Gerhardt}
\address{Ruprecht-Karls-Universit\"at, Institut f\"ur Angewandte Mathematik,
Im Neuenheimer Feld 294, 69120 Heidelberg, Germany}
\email{gerhardt@math.uni-heidelberg.de}
\urladdr{http://www.math.uni-heidelberg.de/studinfo/gerhardt/}
\thanks{This work has been supported by the Deutsche Forschungsgemeinschaft.}

%
\subjclass[2000]{35J60, 53C21, 53C44, 53C50, 58J05}
\keywords{curvature estimates, Weingarten hypersurface, curvature flows}
\date{\today}
%


\begin{abstract}
We prove curvature estimates for general curvature functions. As an application we show the existence of closed, strictly convex hypersurfaces with prescribed  curvature $F$, where the defining cone of $F$ is $\C_+$. $F$ is only assumed to be monotone, symmetric, homogeneous of degree $1$, concave and of class $C^{m,\al}$, $m\ge4$.
\end{abstract}

\maketitle

\tableofcontents

\setcounter{section}{0}
\section{Introduction}
Let $N=N^{n+1}$ be a Riemannian manifold, $\Om\su N$ open, connected and precompact, and $M\su \Om$ a closed connected hypersurface with second fundamental form $h_{ij}$, induced metric $g_{ij}$ and principal curvatures $\ka_i$. $M$ is said to be a Weingarten hypersurface, if, for a given curvature function $F$, its principal curvatures lie in the convex cone $\C\su \R[n]$ in which the curvature function is defined, $M$ is then said to be \tit{admissible}, and satisfies the equation
\begin{equation}\lae{1.1}
\fv FM=f
\end{equation}
where the right-hand side $f$ is a prescribed positive function defined in $\bar\Om$. 

When proving a priori estimates for solutions of \re{1.1} the concavity of $F$ plays a central role. As usual we consider $F$ to be defined in a cone $\C$ as well as on the space of admissible tensors such that 
\begin{equation}
F(h_{ij})=F(\ka_i).
\end{equation}
Notice that curvature functions are always assumed to be symmetric and if $F\in C^{m,\al}(\C)$, $2\le m$, $0<\al<1$, then $F\in C^{m,\al}(\mc S_\C)$, where $\mc S_\C\su T^{0,2}(M)$ is the open set of admissible symmetric tensors with respect to the given metric $g_{ij}$. The result is due to Ball, \cite{ball}, see also \cite[Theorem 2.1.8]{cg:cp}.

The second derivatives of $F$ then satisfy
\begin{equation}
\begin{aligned}
F^{ij,kl}\h_{ij}\h_{kl}=\sum_{i,j}\frac{\pa^2F}{\pa\ka_i\pa\ka_j}\h_{ii}\h_{jj}+\sum_{i\ne j} \frac{F_i-F_j}{\ka_i-\ka_j}(\h_{ij})^2\le 0\;\A\, \h\in \mc S,
\end{aligned}
\end{equation}
where $\mc S\su T^{0,2}(M)$ is the space of symmetric tensors, if $F$ is concave in $\C$, \cf \cite[Lemma 1.1]{cg96}.

However, a mere non-positivity of the right-hand side is in general not sufficient to prove a priori estimates for the $\ka_i$ resulting in the fact that only for special curvature functions for which a stronger estimate was known such a priori estimates could be derived and the problem \re{1.1} solved, if further assumptions are satisfied. 

Sheng et al.\ then realized in \cite{urbas:duke} that the term
\begin{equation}
\sum_{i\ne j} \frac{F_i-F_j}{\ka_i-\ka_j}(\h_{ij})^2
\end{equation}
was all that is needed to obtain the stronger concavity estimates under certain circumstances. Indeed, if the $\ka_i$ are labelled 
\begin{equation}\lae{1.5}
\ka_1\le\cdots\le \ka_n,
\end{equation}
then there holds:
\bl
Let $F$ be concave and monotone, and assume $\ka_1<\ka_n$, then
\begin{equation}\lae{1.8}
\sum_{i\ne j} \frac{F_i-F_j}{\ka_i-\ka_j}(\h_{ij})^2\le \frac2{\ka_n-\ka_1}\sum_{i=1}^n(F_n-F_i)(\h_{ni})^2
\end{equation}
for any symmetric tensor $(\h_{ij})$, where we used coordinates such that $g_{ij}=\de_{ij}$.
\el

\bp
Without loss of generality we may assume that the $\ka_i$ satisfy the strict inequalities
\begin{equation}
\ka_1<\cdots<\ka_n,
\end{equation}
since these points are dense. The concavity of $F$ implies
\begin{equation}
F_1\ge\cdots\ge F_n,
\end{equation}
\cf \cite[Lemma 2]{eh2}, where
\begin{equation}
F_i=\pde F{\ka_i}>0;
\end{equation}
the last inequality is the definition of monotonicity. 
The inequality then follows immediately.
\ep

The right-hand side of inequality \re{1.8} is exactly the quantity that is needed to balance a bad technical term in the a priori estimate for $\ka_n$, at least in Riemannian manifolds, as we shall prove. Unfortunately, this doesn't work in Lorentzian spaces, because of a sign difference in the Gau{\ss} equations.

The assumptions on the curvature function are very simple.
\bas\laas{1.2}
Let $\C\su\R[n]$ be an open, symmetric, convex cone containing $\C_+$ and let $F\in C^{m,\al}(\C)\ii C^0(\bar\C)$, $m\ge 4$, be  symmetric, monotone, homogeneous of degree $1$, and concave such that
\begin{equation}
F>0\qq\text{in}\q\C
\end{equation}
and 
\begin{equation}\lae{1.11}
\fv F{\pa\C}=0.
\end{equation}
\eas

These conditions on the curvature function will suffice. They could have been modified, even relaxed, e.g.,  by  only requiring that $\log F$ is concave, but then the condition 
\begin{equation}
F^{ij}g_{ij}\ge c_0>0,
\end{equation}
which automatically holds, if $F$ is concave and homogeneous of degree  $1$, would have been added, destroying the aesthetic simplicity of \ras{1.2}.

Our estimates apply equally well to solutions of an equation as well as to solutions of curvature flows. Since curvature flows encompass equations, let us state the main estimate for curvature flows.

Let $\Om\su N$ be precompact and connected, and $0<f\in C^{m,\al}(\bar\Om)$. We consider the curvature flow
\begin{equation}\lae{1.13} 
\begin{aligned}
\dot x&=-(\F-\tilde f)\nu\\
x(0)&=x_0,
\end{aligned}
\end{equation}
where $\F$ is  $\F(r)=r$ and $\tilde f=f$, $x_0$ is the embedding of an initial admissible hypersurface $M_0$ of class $C^{m+2,\al}$ such that 
\begin{equation}
\F-\tilde f\ge 0\qq\text{at}\q t=0,
\end{equation}
where of course $\F=\F(F)=F$. We introduce the technical function $\F$ in the present case only to make a comparison with former results, which all use the notation for the more general flows, easier.

We assume that $\bar\Om$ is covered by a Gaussian coordinate system $(x^\al)$, $0\le 1\le n$, such that the metric can be expressed as 
\begin{equation}
d\bar s^2=e^{2\psi}\{(dx^0)^2+\s_{ij}dx^idx^j\}
\end{equation}
and $\bar\Om$ is covered by the image of the cylinder
\begin{equation}
I\times \so
\end{equation}
where $\so$ is a compact Riemannian manifold and $I=x^0(\bar\Om)$,  $x^0$ is a global coordinate defined in $\bar\Om$ and $(x^i)$ are local coordinates of $\so$.

Furthermore we assume that $M_0$ and the other flow hypersurfaces can be written as graphs over $\so$. The flow should exist in a maximal time interval $[0,T^*)$, stay in $\Om$, and uniform $C^1$-estimates should already have been established.

\br
The assumption on the existence of the Gaussian coordinate system and the fact that the hypersurfaces can be written as graphs could be replaced by assuming the existence of a unit vector field $\h\in C^2(T^{0,1}(\bar\Om))$ and of a constant $\theta>0$ such that
\begin{equation}\lae{1.17}
\spd\h\nu\ge 2\theta
\end{equation}
uniformly during the flow, since this assumption would imply uniform $C^1$-estimates, which are the requirement that the induced metric can be estimated accordingly by  controlled metrics from below and above, and because the existence of such a vector field is essential for the curvature estimate.

If the flow hypersurfaces are graphs in a Gaussian coordinate system, then such a vector field is given by
\begin{equation}\lae{1.18}
\h=(\h_\al)=e^{\psi}(1,0,\dots,0)
\end{equation}
and the $C^1$-estimates are tantamount to the validity of inequality \re{1.17}.

In case $N=\R$ and starshaped hypersurfaces one could also use the term
\begin{equation}
\spd x\nu,
\end{equation}
\cf \cite[Lemma 3.5]{cg90}.
\er

Then we shall prove: 
\bt\lat{1.4}
Under the assumptions stated above the principal curvatures $\ka_i$ of the flow hypersurfaces are uniformly bounded from above
\begin{equation}
\ka_i\le c,
\end{equation}
provided there exists a strictly convex function $\chi\in C^2(\bar\Om)$. The constant $c$ only depends on $\abs{f}_{2,\Om}$, $\theta$, $F(1,\dots,1)$, the initial data, and the estimates for $\chi$  and those of the ambient Riemann curvature tensor in $\bar\Om$.

Moreover, the $\ka_i$ will stay in a compact set of $\C$.
\et

As an application of this estimate our former results on the existence of a strictly convex hypersurface $M$ solving the equation \re{1.1}, \cite{cg96,cg97}, which we proved for curvature functions $F$ of class $(K)$, are now valid for curvature functions $F$ satisfying \ras{1.2} with $\C=\C_+$. 

We are even able to solve the existence problem by using a curvature flow which formerly only worked in case that the sectional curvature of the ambient space was non-positive.
\bt\lat{1.5}
Let $F$ satisfy the assumptions above with $\C=\C_+$ and assume that the boundary of $\Om$ has two components
\begin{equation}
\pa\Om=M_1\uud M_2,
\end{equation}
where the $M_i$ are closed, connected strictly convex hypersurfaces of class $C^{m+2,\al}$, $m\ge 4$, which can be written as graphs in a normal Gaussian coordinate system covering $\bar\Om$, and where we assume that the normal of $M_1$ points outside of $\Om$ and that of $M_2$ inside. Let $0<f\in C^{m,\al}(\bar\Om)$,  and assume that $M_1$ is a lower barrier for the pair $(F,f)$ and $M_2$ an upper barrier, then the problem \re{1.1} has a strictly convex solution $M\in C^{m+2,\al}$ provided there exists a strictly convex function $\chi\in C^2(\bar\Om)$. The solution is the limit hypersurface of a converging curvature flow.
\et

\section{Curvature estimates}
Let $M(t)$ be the flow hypersurfaces, then their second fundamental form $h^j_i$ satisfies the evolution equation, \cf \cite[Lemma 2.4.1]{cg:cp}:
\bl\lal{1.4.1}
The mixed tensor $h_i^j$ satisfies the parabolic equation
\begin{equation}\lae{1.4.1}
\begin{aligned}
&\qq\qq\dot h_i^j-\dot\F F^{kl}h_{i;kl}^j=\\[\cma]
&\hp{=}\; \dot\F F^{kl}h_{rk}h_l^rh_i^j-\dot\F F
h_{ri}h^{rj}+ (\F-\tilde f) h_i^kh_k^j\\
&\hp{+}-\tilde f_{\alpha\beta} x_i^\alpha x_k^\beta g^{kj}+ \tilde f_\alpha\n^\alpha h_i^j+\dot\F
F^{kl,rs}h_{kl;i}h_{rs;}^{\hphantom{rs;}j}\\
&\hp{=}+\ddot \F F_i F^j+2\dot \F F^{kl}\riema \alpha\beta\gamma\delta x_m^\alpha x_i ^\beta x_k^\gamma
x_r^\delta h_l^m g^{rj}\\
&\hp{=}-\dot\F F^{kl}\riema \alpha\beta\gamma\delta x_m^\alpha x_k ^\beta x_r^\gamma x_l^\delta
h_i^m g^{rj}-\dot\F F^{kl}\riema \alpha\beta\gamma\delta x_m^\alpha x_k ^\beta x_i^\gamma x_l^\delta h^{mj} \\
&\hp{=}+\dot\F F^{kl}\riema \alpha\beta\gamma\delta\n^\alpha x_k^\beta\n^\gamma x_l^\delta h_i^j-\dot\F F
\riema \alpha\beta\gamma\delta\n^\alpha x_i^\beta\n^\gamma x_m^\delta g^{mj}\\
&\hp{=}+ (\F-\tilde f)\riema \alpha\beta\gamma\delta\n^\alpha x_i^\beta\n^\gamma x_m^\delta g^{mj}\\
&\hp{=}+\dot\F F^{kl}\bar R_{\alpha\beta\gamma\delta;\e}\{\n^\alpha x_k^\beta x_l^\gamma x_i^\delta
x_m^\e g^{mj}+\n^\alpha x_i^\beta x_k^\gamma x_m^\delta x_l^\e g^{mj}\}.
\end{aligned}
\end{equation}
\el

Let $\h$ be the vector field \re{1.18}, or any vector field satisfying \re{1.17}, and set
\begin{equation}
\tilde v=\spd\h\nu,
\end{equation}
then we have:
\bl[Evolution of $\tilde v$]\lal{1.4.4}
The quantity $\tilde v$ satisfies the evolution equation
\begin{equation}\lae{1.4.23}
\begin{aligned}
\dot{\tilde v}-\dot\F F^{ij}\tilde v_{ij}=&\dot\F F^{ij}h_{ik}h_j^k\tilde v
-[(\F-\tilde f)-\dot\F F]\h_{\alpha\beta}\n^\alpha\n^\beta\\
&-2\dot\F F^{ij}h_j^k x_i^\alpha x_k^\beta \h_{\alpha\beta}-\dot\F F^{ij}\h_{\alpha\beta\gamma}x_i^\beta
x_j^\gamma\n^\alpha\\
&-\dot\F F^{ij}\riema \alpha\beta\gamma\delta\n^\alpha x_i^\beta x_k^\gamma x_j^\delta\h_\e x_l^\e g^{kl}\\
&-\tilde f_\beta x_i^\beta x_k^\alpha \h_\alpha g^{ik}.
\end{aligned}
\end{equation}
\el

The derivation is elementary, see the proof of the corresponding lemma in the Lorentzian case \cite[Lemma 2.4.4]{cg:cp}.

Notice that $\tilde v$ is supposed to satisfy \re{1.17}, hence
\begin{equation}
\f=-\log(\tilde v-\theta)
\end{equation}
is well defined and there holds
\begin{equation}\lae{2.5}
\begin{aligned}
\dot \f-\dot\F F^{ij}\f_{ij}=-\{\dot {\tilde v}-\dot\F F^{ij}\tilde v_{ij}\}\frac1{\tilde v-\theta}-\dot\F F^{ij}\f_i\f_j.
\end{aligned}
\end{equation}

Finally, let $\chi$ be the strictly convex function. Its evolution equation is
\begin{equation}\lae{2.6}
\begin{aligned}
\dot\chi-\dot\F F^{ij}\chi_{ij}&=-[(\F-\tilde f)-\dot \F F]\chi_\al\nu^\al -\dot\F F^{ij}\chi_{\al\bet}x^\al_ix^\bet_j\\
&\le -[(\F-\tilde f)-\dot \F F]\chi_\al\nu^\al -c_0 \dot\F F^{ij}g_{ij}
\end{aligned}
\end{equation}
where $c_0>0$ is independent of $t$.

We can now prove \rt{1.4}:
\bp[Proof of \rt{1.4}]
Let $\zeta$ and $w$ be respectively defined by
\begin{align}
\zeta&=\sup\set{{h_{ij}\h^i\h^j}}{{\norm\h=1}},\\
w&=\log\zeta+ \f+\lam \chi,\lae{2.3.14}
\end{align}
where $\lam>0$ is supposed to be large.  We claim that
$w$ is bounded, if $\lam$ is chosen sufficiently large.

Let $0<T<T^*$, and $x_0=x_0(t_0)$, with $ 0<t_0\le T$, be a point in $M(t_0)$ such
that
\begin{equation}
\sup_{M_0}w<\sup\set {\sup_{M(t)} w}{0<t\le T}=w(x_0).
\end{equation}

We then introduce a Riemannian normal coordinate system $(\x^i)$ at $x_0\in
M(t_0)$ such that at $x_0=x(t_0,\x_0)$ we have
\begin{equation}
g_{ij}=\delta_{ij}\q \tup{and}\q \zeta=h_n^n.
\end{equation}

Let $\tilde \h=(\tilde \h^i)$ be the contravariant vector field defined by
\begin{equation}
\tilde \h=(0,\dotsc,0,1),
\end{equation}
and set
\begin{equation}
\tilde \zeta=\frac{h_{ij}\tilde \h^i\tilde \h^j}{g_{ij}\tilde \h^i\tilde \h^j}\raise 2pt
\hbox{.}
\end{equation}

$\tilde\zeta$ is well defined in neighbourhood of $(t_0,\x_0)$.

Now, define $\tilde w$ by replacing $\zeta$ by $\tilde \zeta$ in \re{2.3.14}; then, $\tilde w$
assumes its maximum at $(t_0,\x_0)$. Moreover, at $(t_0,\x_0)$ we have 
\begin{equation}
\dot{\tilde \zeta}=\dot h_n^n,
\end{equation}
and the spatial derivatives do also coincide; in short, at $(t_0,\x_0)$ $\tilde \zeta$
satisfies the same differential equation \re{1.4.1} as $h_n^n$. For the sake of
greater clarity, let us therefore treat $h_n^n$ like a scalar and pretend that $w$
is defined by 
\begin{equation}
w=\log h_n^n+ \f+\lam \chi.
\end{equation} 

From the equations \re{1.4.1}, \re{2.5},  \re{2.6} and \re{1.8}, we infer, by observing the special form of $\F$, i.e., $\F(F)=F$, $\dot\F=1$, $\tilde f=f$ and using the monotonicity and homgeneity of $F$
\begin{equation}
F=F(\ka_i)=F(\tfrac{\ka_1}{\ka_n},\ldots,1)\ka_n\le F(1,\dots,1)\ka_n
\end{equation}
that in $(t_0,\xi_0)$
\begin{equation}\lae{2.15}
\begin{aligned}
0&\le-\tfrac12\dot\F F^{ij}h_{ki}h^k_k\frac\theta{\tilde v-\theta}-f h^n_n+c(\theta)\dot\F F^{ij}g_{ij}+\lam c\\
&\hp{\le}\;-\lam c_0\dot\F F^{ij}g_{ij}-\dot \F F^{ij}\f_i\f_j+\dot\F F^{ij}(\log h^n_n)_i(\log h^n_n)_j\\
&\hp{\le}\;+\frac2{\ka_n-\ka_1}\dot\F \sum_{i=1}^n(F_n-F_i)(h_{ni;}^{\hp{ni;}n})^2 (h^n_n)^{-1}.
\end{aligned}
\end{equation}
Similarly as in \cite[p. 197]{cg:scalar}, we distinguish two cases

\cvm
\tit{Case} $1$.\q Suppose that
\begin{equation}
\abs{\ka_1}\ge \e_1 \ka_n,
\end{equation}
where $\e_1>0$ is small, notice that the principal curvatures are labelled according to \re{1.5}. Then, we infer from \cite[Lemma 8.3]{cg:scalar}
\begin{equation}
F^{ij}h_{ki}h^k_j\ge \tfrac1n F^{ij}g_{ij}\e_1^2\ka_n^2,
\end{equation}
and 
\begin{equation}
F^{ij}g_{ij}\ge F(1,\ldots,1),
\end{equation}
for a proof see e.e., \cite[Lemma 2.2.19]{cg:cp}.

Since $Dw=0$,
\begin{equation}
D\log h^n_n=-D\f-\lam D\chi,
\end{equation}
we obtain
\begin{equation}
\dot\F F^{ij}(\log h^n_n)_i(\log h^n_n)_j=\dot \F F^{ij}\f_i\f_j+2\lam \dot\F F^{ij}\f_i\chi_j+\lam^2\dot\F F^{ij}\chi_i\chi_j,
\end{equation}
where
\begin{equation}
\abs{\f_i}\le c\abs{\ka_i}+c,
\end{equation}
as one easily checks. 

Hence, we conclude that $\ka_n$ is a priori bounded in this case.

\cvm
\tit{Case} $2$.\q Suppose that
\begin{equation}\lae{2.22}
\ka_1\ge -\e_1\ka_n,
\end{equation}
then, the last term in inequality \re{2.15} is estimated from above by
\begin{equation}
\begin{aligned}
&\frac2{1+\e_1}\dot\F \sum_{i=1}^n(F_n-F_i)(h_{ni;}^{\hp{ni;}n})^2 (h^n_n)^{-2}&\le  \\
&\frac2{1+2\e_1}\dot\F \sum_{i=1}^n(F_n-F_i)(h_{nn;}^{\hp{nn;}i})^2 (h^n_n)^{-2}\\
&\qq +c(\e_1)\dot\F \sum_{i=1}^{n-1}(F_i-F_n)\ka_n^{-2}
\end{aligned}
\end{equation} 
where we used the Codazzi equation. The last sum can be easily balanced.

The terms in \re{2.15} containing the derivative of $h^n_n$ can therefore be estimated from above by
\begin{equation}
\begin{aligned}
&-\frac{1-2\e_1}{1+2\e_1}\dot\F \sum_{i=1}^nF_i(h_{nn;}^{\hp{nn;}i})^2 (h^n_n)^{-2}\\
&+\frac2{1+2\e_1}\dot\F F_n\sum_{i=1}^n(h_{nn;}^{\hp{nn;}i})^2 (h^n_n)^{-2}\\
&\le\dot\F F_n\sum_{i=1}^n(h_{nn;}^{\hp{nn;}i})^2 (h^n_n)^{-2}\\
&=\dot\F F_n \norm{D\f+\lam D\chi}^2\\
&=\dot\F F_n\{\norm{D\f}^2+\lam^2\norm{D\chi}^2+2\lam \spd{D\f}{D\chi}\}.
\end{aligned}
\end{equation}

Hence we finally deduce
\begin{equation}
\begin{aligned}
0\le -\dot\F \tfrac12 F_n\ka_n^2\frac\theta{\tilde v-\theta}&+c\lam^2\dot\F F_n(1+\ka_n)-f\ka_n+\lam c \\
&+(c(\theta)-\lam c_0)\dot\F F^{ij}g_{ij}
\end{aligned}
\end{equation}
Thus, we obtain an a priori estimate 
\begin{equation}\lae{2.28}
\ka_n\le \const,
\end{equation}
if $\lam$ is chosen large enough. Notice that $\e_1$ is only subject to the requirement $0<\e_1<\frac12$.
\ep

\br
Since  the initial condition $F\ge f$ is preserved under the flow, a simple application of the maximum principle, \cf \cite[Lemma 5.2]{cg96}, we conclude that the principal curvatures of the flow hypersurfaces stay in a compact subset of $\C$.
\er

\br\lar{2.3}
These a priori estimates are of course also valid, if $M$ is a stationary solution.
\er

\section{Proof of \rt{1.5}}

We  consider the curvature flow \re{1.13} with initial hypersurface $M_0=M_2$. The flow will exist in a maximal time interval $[0,T^*)$ and will stay in $\bar\Om$. We shall also assume that $M_2$ is not already a solution of the problem for otherwise the flow will be stationary from the beginning.

Furthermore, the flow hypersurfaces can be written as graphs
\begin{equation}
M(t)=\graph u(t,\cdot)
\end{equation}
over $\so$, since the initial hypersurface has this property and all flow hypersurfaces are supposed to be convex, i.e., uniform $C^1$-estimates are guaranteed, \cf \cite{cg96}.

The curvature estimates from \rt{1.4} ensure that the curvature operator is uniformly elliptic, and in view of  well-known regularity results we then conclude that the flow exists for all time and converges in $C^{m+2,\bet}(\so)$ for some $0<\bet\le\al$ to a limit hypersurface $M$, that will be a stationary solution, \cf \cite[Section 6]{cg:survey}.


\providecommand{\bysame}{\leavevmode\hbox to3em{\hrulefill}\thinspace}
\providecommand{\href}[2]{#2}



\end{document}